\documentstyle[12pt]{article} 
\newcommand{\r}{{\iota}}

\newcommand{\Uq}{U_q(g)} 
\newcommand{\Ur}{U_q(g)} 
\newcommand{\E}{\cal E}   
\newcommand{\lam}{\lambda} 
\newcommand{\Lam}{\Lambda} 
\newcommand{\x}{\otimes}

\newcommand{\ba}{\begin{eqnarray}} 
\newcommand{\na}{\end{eqnarray}} 
\newcommand{\ban}{\begin{eqnarray*}} 
\newcommand{\nan}{\end{eqnarray*}}

\newcommand{\K}{K_{2\rho}}

\begin{document} 
\title{\small{ TOPOLOGICAL INVARIANTS FOR LENS SPACES AND\\ 
EXCEPTIONAL  QUANTUM GROUPS}  }   
\author{ R. B. Zhang\\   
\small  Department of Pure Mathematics\\ 
\small  University of Adelaide\\ 
\small  Adelaide, Australia} 
\date{November 1995}
\maketitle  

\vspace{3cm}  
The  Reshetikhin - Turaev invariants arising from 
the quantum groups associated with the exceptional Lie algebras  
$G_2$, $F_4$ and $E_8$ at odd roots of unity are constructed and 
explicitly computed for all the lens spaces.  
  
\vspace{3cm}

\section{Introduction}  

A pressing problem in the field of `quantum topology' 
\cite{Kauffman}\cite{Turaev}  
is to understand the topological information embodied in the 
quantum invariants of $3$-manifolds\cite{Witten} - 
\cite{Kohno}  constructed in recent years,
and to use the information to settle geometric questions.  
A direct way to tackle the problem is to compute these invariants 
for $3$-manifolds of interest, then try to 
extract topological information from the explicit results.  
This approach was followed by a number of authors, and 
extensive studies were carried out for the $su(2)$ invariants 
of Witten - Reshetikhin - Turaev, though only for classes of 
$3$-manifolds of a very simple nature \cite{Singer}-\cite{Murakami}.    

Even in the case of the smallest Lie algebra $su(2)$, 
the computations of such invariants is a rather difficult task,  
as it involves the evaluation of link invariants 
arising from all the irreducible representations 
of the $sl(2)$ quantum group at roots of unity.   
For bigger Lie algebras, the level of difficulty increases 
considerably.  Nevertheless, we succeeded in computing the invariants 
for the lens spaces using all the classical Lie algebras \cite{Wenzl}
in an earlier publication \cite{Carey}. 
It is the aim of this letter to extend the 
results of \cite{Carey} to the exceptional Lie algebras.  
We will be mainly concerned with $G_2$, $F_4$ and $E_8$. 
The reason we work with these Lie algebras is that 
the fundamental groups of their root systems are all trivial. 
This fact greatly simplifies the computations. 
The main result of this letter is an 
explicit formula for the Reshetikhin-Turaev invariants arising 
from the quantum groups associated with the above mentioned 
Lie algebras at odd roots of unity  for the lens spaces. 
We also present a simple construction of these  
invariants, which enables us to carry out
the computations relatively easily.

Recall that the  Reshetikhin-Turaev construction of 
$3$-manifold invariants using quantum groups is usually 
formulated in terms of modular Hopf algebras \cite{Reshetikhin}
\cite{Turaev}.  A modular Hopf 
algebra is a ribbon Hopf algebra admitting a distinguished set of 
representations, which satisfy some very rigid conditions. 
It is in general a very difficult problem to show that a Hopf algebra 
satisfies the requirements of a modular Hopf algebra.  In \cite{Zhang}, 
we proposed a slightly different construction,   
which did not depend on the notion of modular Hopf algebras. 
Instead, it relied on the analysis of eigenvalues of certain central 
elements of the quantum groups under consideration.   
The construction worked rather  simply for the quantum supergroups
$U_q(osp(1|2))$ and $U_q(gl(1|2))$, and the classical 
series of quantum groups.  As we will see below, it can be 
easily implemented for the exceptional quantum groups as well.

The organization of this letter is as follows. In section 2 we present 
some results on the representation theory and the central algebra of 
the quantum groups.  In section 3 we construct the topological invariants 
then compute them for the lens spaces.

\section{Exceptional Quantum Groups } 
We will consider only the exceptional Lie algebras $G_2$, $F_4$ and $E_8$.  
Let $g$ be any of these Lie algebras, and denote the quantum group 
associated with $g$ by $\Uq$. 
Denote by $\Phi^+$ the set of the positive roots of $g$ relative to 
a base  $\Pi=\{\alpha_1, ..., \alpha_r\}$. 
Define $H^*=\bigoplus_{i=1}^\r{\bf C}\alpha_i$,  
$\E=\bigoplus_{i=1}^r{\bf R}\alpha_i$.    Let $(\ , \ ): \E\times \E 
\rightarrow {\bf R}$ be an inner product of $\E$  such that the Cartan 
matrix $A$ of $g$ is given by 
\ban 
A=\left( a_{i j}\right)_{i j=1}^r, & a_{i j}={ {2(\alpha_i, \alpha_j)} 
\over{(\alpha_i, \alpha_i)}}.  
\nan  
We will normalize this inner product such that,   
for the highest root $\theta$ of $g$,  
\ban 
(\theta, \ \theta)&=&\left\{\begin{array}{l l} 
                           6, & for \ G_2, \\ 
                           4, &for \  F_4, \\ 
                           2, &for \  E_8. 
                       \end{array}\right. 
\nan  
We assume that $q$ is a primitive $N$-th root of unity, 
where $N$ is a positive odd integer, which is greater than 
the dual Coxeter number of $g$, and is not divisible by $3$ in the 
case of $G_2$.

The quantum group \cite{Drinfeld}  $\Uq$ is defined to be the 
unital associative algebra generated by $\{ k_i, k_i^{-1}, e_i, f_i\ | \ 
i=1, ..., r\}$, $r$ being the rank of $g$,  with the following relations 
\ba 
k_i k_j = k_j k_i, & k_i k_i^{-1} =1, \nonumber\\ 
k_i e_j k_i^{-1} = q_i^2 e_j, & 
k_i f_j k_i^{-1} = q_i^{-2} f_j,   
\na 
\ba 
{[}e_i, \ f_j]=\delta_{i j}{{ k_i - k_i^{-1}}
\over{q_i -q_i^{-1}}}, &&\nonumber\\   
\sum_{t=0}^{1-a_{i j}} (-1)^t 
\left[ \begin{array}{r r r} 
        1&-&a_{i j} \\ 
           &t&    
        \end{array}\right]_{q_i}        
(e_i)^t e_j (e_i)^{1-a_{i j}-t}&=& 0 , \ \ \ \ i\ne j,  \nonumber\\ 
\sum_{t=0}^{1-a_{i j}} (-1)^t 
\left[ \begin{array}{r r r}     
        1&-&a_{i j} \\ 
           &t&     
        \end{array}\right]_{q_i}     
(f_i)^t f_j (f_i)^{1-a_{i j}-t}&=& 0 , \ \ \ \ i\ne j,   
\na 
where $\left[ \begin{array}{r}
       s\\ t  
        \end{array}\right]_q $  
is the Gauss polynomial, and  $q_i=q^{(\alpha_i, \ \alpha_i)/2}$.

It is well known that the quantum group $\Uq$ has the structure of 
a Hopf algebra. We will take the following co-multiplication 
\ban 
\Delta(k_i^{\pm 1})&=&k_i^{\pm 1}\otimes k_i^{\pm 1},\\  
\Delta(e_i)&=&e_i\otimes k_i + 1\otimes e_i, \\    
\Delta(f_i)&=&f_i\otimes 1 + k_i^{-1}\otimes f_i.     
\nan 
A co-unit and an antipode can also be defined, but we will not spell them
out explicitly here. 

As it stands, $\Uq$ does not admit a universal $R$ - matrix, 
but an appropriate Hopf quotient of it does.  
To explain this quotient, we follow \cite{Lusztig} to 
define  the  elements of $\Uq$:   
$e^{(n)}_\alpha,  f^{(n)}_\alpha$,  $\alpha\in\Phi^+$, $n=0, 1, ...$,  
where $e^{(1)}_\alpha$ and  $f^{(1)}_\alpha$ have properties 
similar to that of the root spaces of the Lie algebra $g$.  
In particular, if $\alpha =\alpha_i$ for a given $i$, 
then $e^{(1)}_\alpha = e_i$, and $f^{(1)}_\alpha=f_i$.  
Then the  elements $e^{(N)}_\alpha,  f^{(N)}_\alpha$, 
$k_i^N - 1$,  $\alpha\in\Phi^+$, $i=1, 2, ..., r$, 
generate a Hopf ideal $J$ of $\Uq$.  
The quotient algebra $\Uq/J$ is again Hopf, and is known to be 
quasi triangular, that is,  it admits a universal $R$-matrix. 
By an abuse of notation, we still denote this quotient algebra 
by $\Uq$, noting that it 
 is clearly finite dimensional over the complex field.

Express the universal $R$ matrix of $\Uq$
as $R=\sum_t  a_t\x b_t$, and define
\ban
u&=&\sum_t S(b_t)a_t.
\nan
Set 
\ban
\K&=& \prod_{i=1}^r (k_i)^{2\rho_i}
\nan
where $\rho_i\in{\bf Z}_+$ are specified by 
$\sum_{i=1}^r\rho_i\alpha_i=\sum_{\alpha\in \Phi^+} \alpha/2=\rho$. 
Then 
\ban
v&=&u\K^{-1},  
\nan
belongs to the center of $\Uq$, and satisfies
\ban
\Delta(v)&=& (v\otimes v)(R^T R)^{-1}.
\nan
 
Let $V$ be a finite dimensional $\Uq$ module.
We denote the corresponding representation by $\pi$.
Define
\ba
C_V&=&tr_V[(\pi\x id)(\K^{-1}\x 1)R^T R],
\label{eq:C}
\na
where $tr_V$ represents the trace taken over $V$.
Then $C_V$ is also central.\\

Since the quantum  group  $\Uq$ is finite dimensional over $\bf C$, 
all its irreducible representations are finite dimensional as well, 
following the textbook result that every irreducible left 
module over an associative algebra is isomorphic 
to the quotient of the algebra itself by a maximal left ideal. 
It can also be shown that every irreducible $\Uq$ 
module admits a unique (up to scalar multiples) highest 
and lowest weight vector.  The highest weight vector $v_+$ of  
an irreducible $\Uq$ module $V$ is defined by  
\ban  
k_i v_+ &=& q_i^{l_i} v_+, \\ 
e_i v_+&=&0,     \ \ \ \ i=1, 2, ..., r,     
\nan  
where $l_i\in\{0, 1, ..., N-1\}$ as required by $k_i^N=1$.

Let $X$ denote the weight lattice of $g$.  For the Lie algebras under 
consideration, $X$ coincides with the root lattice of  $g$. 
Now from the above discussion we easily see that 
each irreducible $\Uq$ module is uniquely characterized 
by an element of the set $X_N:=X/ N X$. We denote the canonical 
projection $X\rightarrow X_N$ by $p$.

Let 
\ban 
F&=&\{ x\in X | \ 0<{{2(\lam + \rho, \alpha)}
\over{(\alpha, \alpha)}}< N,  \ \forall \alpha\in\Phi^+\}, \\  
\overline{F}&=&\{ x\in X | \ 0\le{{2(\lam + \rho, \alpha)}
\over{(\alpha, \alpha)}}\le  N,  \ \forall \alpha\in\Phi^+\}, \\ 
\Lam^+_N&=&p(F), \\  
\overline{\Lam^+_N}&=&p(\overline{F}). 
\nan 
Define the set 
\ban 
{\cal V}(\Lam^+_N)& = &\{V(\lam)\  |\  \lam\in \Lam^+_N\}, 
\nan 
of all the irreducible $\Uq$ modules with highest weights
belonging to $\Lam^+_N$.  
A recent result of Andersen and Paradowski\cite{Andersen} asserts  that  
{\em the tensor product of any finite number of irreducible $\Uq$ modules
$V(\lam^t)$ with highest weights $\lam^t\in \Lam^+_N$ can be decomposed into
\ban
\bigotimes_t V(\lam^t) &=&\bigoplus_{\lam\in \Lam^+_N}
V(\lam)^{\oplus m(\lam)}\oplus \cal N,
\nan
where $m(\lam)$ is the multiplicity of $V(\lam)$ appearing in the
tensor product.  The module $\cal N$ is a direct sum of indecomposable
$\Uq$ modules, which has the property that for any module homomorphism
$f: \cal N\rightarrow \cal N$, the quantum trace of $f$ vanishes
identically, i.e.,
\ban tr_{\cal N}(K_{2\rho}f)&=&0. \nan }

A further useful property of the set $\Lam^+_N$ is the following. 
For any $\lam \in \Lam^+_N$,  we define $\lam^*$ by taking any 
representative of $\tilde\lam\in X$, and setting 
$\lam^* = p(-\tau(\tilde\lam))$, where $\tau$ is the maximal 
element of the Weyl group of $g$.  Clearly $\lam^*$ 
is independent of the representative $\tilde\lam$ chosen.   
It is easy to prove that $\lam^*\in \Lam^+_N$, and 
the dual module $V(\lam)^*$ of $V(\lam)$ is isomorphic to 
$V(\lam^*)$.   \\

Define a function 
$q^{(\ ,\ )}:  X_N\times X_N\rightarrow {\bf C}$ by requiring that its 
pullback by the canonical projection $X\rightarrow X_N$ be given 
by the function $X\times X \rightarrow {\bf C}$,  
$\{ \tilde\lam,  \tilde\mu\} \mapsto q^{(\tilde\lam, \tilde\mu)}$. 
Acting on any $V(\lam)\in \cal V(\Lam^+_N)$, 
the central elements $v$ and $C_\mu:=C_{V(\mu)}$, $\mu\in\Lam^+_N$, 
take the following eigenvalues respectively 
\ban 
\chi_\lam(v) &=& q^{-(\lam + 2\rho, \ \lam)}, \\ 
\chi_\lam(C_\mu) &=&S_{\lam \mu}/Q(\mu),    
\nan  
where 
\ban
Q(\mu)&=&\sum_{\sigma\in{\cal W}} \det\sigma q^{2(\sigma(\rho), \ \mu+\rho)},\\
S_{\lam \mu} &=&\sum_{\sigma\in{\cal W}}   \det\sigma \
q^{2(\sigma(\lam+\rho), \ \mu+\rho)},
\nan
In these equations, $\rho$ should be interpreted as $p(\rho)$. 
The $\cal W$ represents the Weyl group of the Lie algebra $g$, 
which  acts on $X_N$ by  
\ban  
\sigma(\tilde\mu + N X) = \sigma(\tilde\mu) + N X, & \tilde\mu\in X, 
& \sigma\in\cal W. 
\nan  
Set 
\ba
d_\lam &=& \Omega Q(\lam), \ \ \ \   \lam\in\Lam^+_N,  \label{solution} \\
\Omega&=&(-1)^{|\Phi^+|}
   { {q^{3(\rho, \rho)}}\over{G_1(q; g)}},  \nonumber
\na
where $G_1(q; g)$ is the $k=1$ case of 
\ban 
G_k(q; g)&=&\sum_{\lam\in X_N} \left(q^{(\lam, \ \lam)}\right)^k. 
\nan 
It is clear that $d_\lam=d_{\lam^*}$. 

Define 
\ba
\delta&=&v-\sum_{\lambda\in\Lam^+_N} d_\lambda
\chi_\lambda(v^{-1})C_\lambda.  \label{delta}
\na
Then 
{\em $\delta$ acts as the zero map on any element of 
$\cal V(\Lam^+_N)$, i.e., in any irreducible $\Uq$ module 
with highest weight belonging to $\Lam^+_N$. } 

To understand the above assertion, observe that 
$\overline{\Lam_N^+}$ furnishes a fundamental domain for the action 
of  $\cal W$ on $X_N$, a fact which can be easily proven by 
studying the action of the affine Weyl group ${\cal W}_N$ on $X$.  
Therefore, for any $\lam$,  
$\mu$ $\in{\Lam_N^+}$  and $\sigma$, $\omega$ $\in{\cal W}$, 
\ban 
\sigma(\lam+\rho)-\rho=\omega(\mu+\rho)-\rho
& iff & \sigma=\omega,\ \  \ \ \lam=\mu.
\nan 
Also note that $S_{\lam \nu}=0$ if
$\lam\in{\overline{\Lam_N^+} - \Lam_N^+}$. 
These properties immediately lead to 
\ban
\sum_{\lam\in\Lam^+_N} d_\lam q^{(\lam+2\rho, \lam)} S_{\lam \mu} 
&=& \sum_{\nu\in X_N} x_\nu
q^{(\nu+2\rho, \nu)} S_{\nu  \mu},\\ 
x_\lam &=& { {q^{3(\rho, \rho)-2(\lam+\rho, \rho)}}\over{G_1(q; g)}}.   
\nan 
Now our assertion follows directly from    
\ban  
\sum_{\nu\in X_N} x_\nu q^{(\nu+2\rho, \nu)} S_{\nu  \mu} 
&=& Q(\mu) q^{-(\mu+2\rho, \mu)}.   
\nan

By using the Andersen-Paradowski result and the properties of 
the central element $\delta$ we can show that 
{\em if $Y$ is a $\Uq$ module which is the tensor product 
of a finite number of elements(not necessarily distinct) of 
${\cal V}(\Lam^+_N)$, then for any module homomorphism 
$f:  Y\rightarrow Y$,  
\ba 
tr_Y(\K\delta f) &=&0, \label{vanish}  
\na  
where $tr_Y$ represents the trace taken over $Y$. }

Define  
\ba 
z&=& \sum_{\lam\in\Lam^+_N} d_\lambda q^{-(\lam+2\rho, \lam)} D_q(\lam), 
\label{z}  
\na 
where $D_q(\lam)$ denotes  the quantum dimension of the irreducible 
$\Uq$ module $V(\lam)$ $\in$ ${\cal V}(\Lam^+_N)$, which is given by 
$D_q(\lam)=Q(\lam)/Q(0)$.  Then we have the following useful formulae 
\ban 
z&=&(-1)^{|\Phi^+|} q^{6(\rho, \ \rho)} 
{{G_{N-1}(q; g)}\over{G_1 (q; g)}} ; \\ 
|z|&=& 1.   
\nan  

\section{$3$ - Manifold Invariants} 
In order to construct  $3$-manifold invariants using the quantum 
group $\Uq$, we need to employ the Reshetikhin-Turaev $F$ 
functor from the category of
coloured ribbon graphs to the category of finite dimensional 
representations of this quantum group. 
A detailed discussion of coloured ribbon graphs and this   
functor can be found in \cite{Turaev}; and we will not repeat it here.  
Instead, we consider as examples the ribbon $(k, \ k)$ 
graphs depicted in Figure 1 to provide some concrete intuition
 about the functor. 

\vspace{1 cm}
\begin{center}
Figure  1.     Figures are available on request. 
\end{center}
\vspace{1 cm}

\noindent 
We colour the ribbons of both Figure $1.a$ and $1.b$ by 
$\{V(\lam_1), V(\lam_2), ..., V(\lam_k)\}$,  
while we colour the annulus of Figure 1.a by  
the irreducible  module $V(\mu)$.  
We denote the resultant coloured ribbon graphs by 
$\phi^{(k)}_\mu$ and $\zeta^{(k)}$ respectively. 
It is straightforward to obtain,  
\ba
F(\phi^{(k)}_\mu)&=&\chi_j(v^{-1})\Delta^{(k-1)}(C_\mu),\nonumber\\  
F(\zeta^{(k)})&=&\Delta^{(k-1)}(v),   \label{Kirby}
\na 
which map 
$V(\lam_1)\otimes V(\lam_2)\otimes...\otimes V(\lam_k)$ to itself.  

The Reshetikhin-Turaev construction makes use of two fundamental
theorems in $3$-manifold theory, due to Lickorish \cite{Lickorish}
and Wallace, and Kirby \cite{Kirby} and Craggs respectively.
The Lickorish-Wallace theorem states that
each framed link in $S^3$ determines a closed, orientable
$3$-manifold, and every such $3$-manifold is obtainable
by surgery along a framed link in $S^3$.
The disadvantage of this description of $3$- manifolds is that
different framed links may yield homeomorphic $3$-manifolds
upon surgery.  This problem was resolved by Kirby \cite{Kirby} and Craggs,
and Fenn and Rourke \cite{Fenn}.
These authors proved that orientation preserving homeomorphism classes
of  closed, orientable $3$-manifolds correspond bijectively
to equivalence classes of framed links in $S^3$, where the equivalence
relation is generated by the Kirby moves.
The essential  idea of \cite{Reshetikhin} is to make appropriate
combinations of isotopy invariants of a framed link embedded in $S^3$,
such that they will be intact under the Kirby moves, and thus qualify
as topological invariants of the $3$-manifold obtained
by surgery along this link.

Let  $M_L$ be a closed, oriented $3$ - manifold, which is given a   
surgery description, namely, represented  
by surgery along a framed link $L$  embedded in $S^3$. 
Assume that the framed link $L$(in blackboard framing) consists 
of $m$ components $L_i, \ i=1, ..., m$. 
It gives rise to a unique ribbon graph by 
extending each component $L_i$ to an annulus, which has 
$L_i$ itself and an $L_i'$ as its edges, where $L_i'$ is a 
parallel copy of $L_i$ such that the linking number between the two 
is equal to the framing number of the latter. We denote this 
ribbon graph by $\Gamma(L)$.

We colour $\Gamma(L)$ by associating with each component 
$L_i$ with a $V(\lam_i)\in{\cal V}(\Lam^+_N)$.  Set     
$c=\{\lam^{(1)} , \lam^{(2)}, ..., \lam^{(m)}\}$, 
where some $\lambda^{(j)}$'s may be equal, and denote by 
${\cal C}(L)$ the set of all the distinct $c$'s. 
The ribbon graph coloured by modules associated with $c$ will be 
denoted by $\Gamma_c(L)$.  The Reshetikhin-Turaev functor applied 
to $\Gamma_c(L)$  yields $F(\Gamma_c)$, which is   
a homomorphism of the trivial $\Uq$ module to itself, i.e., 
a complex number. 

Define  
\ba
\Sigma(L)&=&\sum_{c\in {\cal C}(L)} \Pi_{i=1}^{m}d_{\lam^{(i)}}
\ F(\Gamma_c(L)).    \label{eq:Sigma}
\na
Since   $d_\lam =d_{\lam^*}$, $\forall \lam\in \Lam^+_N$,   
{\em $\Sigma(L)$ is independent of the orientation chosen for $L$}. 
It follows from (\ref{vanish}) and (\ref{Kirby}) that 
{\em $\Sigma(L)$ is invariant under the positive Kirby moves depicted in 
Figure 2. a and Figure 2. b. }\\    

\vspace{1cm} 
\begin{center} 
Figure 2.    
\end{center} 
\vspace{1cm}

On the other hand, $\Sigma$ is not invariant under the Kirby $(-)$ 
moves given in Figure 2.c and 2.d.  In particular, 
if $L'$ is the framed link obtained by applying once the special 
Kirby $(-)$ move Figure 2.c, namely, adding a framing $-1$ 
unknot,  to the framed link $L$, then  
\ban 
\Sigma(L')&=&z\Sigma(L).  
\nan 
Since $|z|=1$, the norm of $\Sigma$ remains 
intact under both the Kirby $(+)$ moves and the special Kirby 
$(-)$ move.  In view of the fact that these moves together generate 
the entire Kirby calculus, we conclude that 
$|\Sigma(L)|$ is a topological invariant of $M_L$.  

Let $A_L=(a_{ij})_{m\times m}$ be  
the intersection form on the second homology group 
$H_2(W_L, \ {\bf Z})$ of $W_L$, where $W_L$ is the $4$ - manifold 
bounded by  $M_L$.  Denote by $sign(A_L)$ the number of nonpositive 
eigenvalues of $A_L$. 
Under the special Kirby $(-)$ move, $sign(A_L)=sign(A_{L'})-1$, 
while the positive Kirby moves leave $sign(A_L)$ unchanged.     
Therefore $z^{-sign(A_L)}\Sigma(L)$ is invariant under 
the positive Kirby moves and the special negative Kirby move. 
Hence,   
{\em the  following  quantities 
\ba
\bigtriangledown(M_L)&=&|\Sigma(L)|^2, \label{norm} \\   
{\cal F}(M_L)&=&z^{-sign(A_L)}\Sigma(L).  \label{eq:F}
\na
are topological invariant of the $3$ - manifold $M_L$. } 
Note that $\bigtriangledown(M_L)$ should be closely related to the 
Turaev-Viro invariant \cite{Viro}.\\   

Now we compute these invariants for the lens spaces $L(m, n)$,  
where $m, n\in {\bf Z}$ are  co-prime. 
We  assume that $0<n<m$.  This exhausts all the possible lens spaces 
apart from $S^3$ and $S^2\times S^1$. These two degenerate cases 
can be easily treated; we have 
\ban 
{\cal F}(S^3)&=&1,\\ 
{\cal F}(S^2\times S^1)&=& { {1}\over {\Omega Q(0)} }.   
\nan

For each $L(m, n)$, there always exists a unique set of integers 
$\{ a_1, ..., a_s\}$ for some $s$ with $2\le a_i\in{\bf Z}$, 
such that the ratio $m/n$ can be expressed as a  continued fraction 
\ban 
{{m}\over{n}}&=& a_1 - {{1}\over{a_2 - {{1}\over{a_3 - 
{{1}\over{...... -  {{1}\over{a_{s-1} -  {{1}\over{a_s} }  } } }} } } } }.   
\nan  
Then the manifold can be obtained by surgery along the framed link, 
Figure 3,  \\  

\vspace{1cm} 
\begin{center}  Figure 3.   \end{center}  
\vspace{1cm} 

\noindent 
with framing numbers of the respective components of the link 
being  $a_1$, $a_2$, ..., $a_s$. 

In order to compute the quantum invariants for $L(m, n)$, we define,  
for $\lam\in X_N$, 
\ban 
\varpi_\lam&=&\left\{ \begin{array}{l l}
                 d_\lam/Q(\lam), & \lam\in\Lam^+_N, \\
                  0, & \lam\not\in\Lam^+_N.
                \end{array}\right.
\nan 
It can be proved that 
\ba
\sum_{\sigma\in{\cal W}} \varpi_{\sigma(\lam+\rho)-\rho}
&=&\left\{ \begin{array}{l l}
          \Omega, & {{2(\lam+\rho, \  \alpha)}\over{ (\alpha, \  \alpha)}}
           \not\equiv 0 (mod N), \ \  \forall \alpha\in\Phi^+, \\
           0, & otherwise.
          \end{array}\right. \label{varpi}
\na 
 
Consider a ribbon graph of the form Figure 3 with $k+1$ components.  
We temporarily assume that the all framing numbers are arbitrary integers, 
except that $a_{k+1}=0$. Cutting the left most component open results 
in a new ribbon graph, which has $k$ annuli and one ribbon. 
Order the annuli from right to left, and colour the $i$ - th annulus 
by the irreducible $\Ur$ module $V(\mu_i)$, with $\mu_i\in\Lam^+_N$.    
Colour the ribbon by $V(\lam)$, $\lam\in\Lam^+_N$.
We denote the resultant coloured ribbon graph by 
$\Gamma_\lam(a_1, ..., a_k;  \mu_1, ..., \mu_k)$, and define 
\ban 
h_\lam^{(k)}(a_1, ..., a_k)&=&Q(\lam) \sum_{\mu_1, ..., \mu_k\in\Lam^+_N}
\prod_{i=1}^k d_{\mu_i} F(\Gamma_\lam(a_1, ..., a_k;  \mu_1, ..., \mu_k)). 
\nan

It is easy to establish the following recursive formula 
\ban
h_\lam^{(k)}(a_1, ..., a_k)&=&
\sum_{\mu\in{X_N}} \varpi_{\mu}
\sum_{\sigma\in{\cal W}} \det\sigma  q^{a_k(\mu+2\rho, \ \mu)}
q^{2(\sigma(\mu+\rho) , \ \lam+\rho)}h_\mu^{(k-1)}(a_1, ..., a_{k-1}),
\nan
with  
\ban 
h_\lam^{(1)}(a_1)&=&
{\Omega}\sum_{\mu_1\in{X_N}} Q(\mu_1) q^{a_1(\mu_1+2\rho,
  \ \mu_1) + 2(\mu_1+\rho,  \ \lam+\rho)}. 
\nan  
Observe that the  
$h_\lam^{(t)}(a_1, ..., a_t)$ are well defined for all 
$\lam\in X_N$, and  more importantly,  
\ban
h_\lam^{(t)}(a_1, ..., a_t)&=& 
\det\sigma\  h_{\sigma(\lam+\rho)-\rho}^{(t)}(a_1, ..., a_t),
\ \ \forall \sigma\in{\cal W}. 
\nan  
Using this property and equation (\ref{varpi}) we obtain  
\ban 
h_\lam^{(k)}(a_1, ..., a_k)&=&
\Omega \sum_{\mu\in{X_N}} q^{a_k(\mu+2\rho, \ \mu)
+2(\mu+\rho , \ \lam+\rho)} h_\mu^{(k-1)}(a_1, ..., a_{k-1}).  
\nan 
Since 
\ban 
{\cal F}(L(m, n))&=& {{1}\over{Q(0)}} h_0^{(s)}(a_1, ..., a_s), 
\nan 
we immediately arrive at 
\ba 
{\cal F}(L(m, n))&=& 
{{\Omega^s}\over{Q(0)}}\sum_{\mu_1, ..., \mu_s\in X_N}
Q(\mu_1) q^{\sum_{i=1}^s \left[a_i (\mu_i + 2\rho, \mu_i) + 2(\mu_i +\rho,
\mu_{i+1}+\rho)\right]}, \label{H}       
\na  
where $\mu_{s+1}=0$.  Also, 
\ba 
\bigtriangledown(L(m, n)))&=& |{\cal F}(L(m, n))|^2. 
\na

\end{document}